\documentclass[conf]{new-aiaa}
\usepackage[utf8]{inputenc}

\usepackage{todonotes}
\usepackage{graphicx}
\usepackage{amsmath}
\usepackage{ulem}
\usepackage[version=4]{mhchem}
\usepackage{siunitx}
\usepackage{longtable,tabularx}
\usepackage{svg}
\usepackage{subcaption}
\title{Energy-Preserving Reduced Operator Inference \\ for Efficient Design and Control}

\author{Tomoki Koike\footnote{PhD Candidate, School of Aerospace Engineering, AIAA Student Member, tkoike3@gatech.edu} and Elizabeth Qian\footnote{Assistant Professor, Schools of Aerospace Engineering \& Computational Science and Engineering, AIAA Member, eqian@gatech.edu}}
\affil{Georgia Institute of Technology, Atlanta, Georgia, 30332}

\begin{document}

\maketitle

\begin{abstract}
Many-query computations, in which a computational model for an engineering system must be evaluated many times, are crucial in design and control. For systems governed by partial differential equations (PDEs), typical high-fidelity numerical models are high-dimensional and too computationally expensive for the many-query setting. 
Thus, efficient surrogate models are required to enable low-cost computations in design and control. This work presents a physics-preserving reduced model learning approach that targets PDEs whose quadratic operators preserve energy, such as those arising in governing equations in many fluids problems. The approach is based on the Operator Inference method, which fits reduced model operators to state snapshot and time derivative data in a least-squares sense. However, Operator Inference does not generally learn a reduced quadratic operator with the energy-preserving property of the original PDE. Thus, we propose a new energy-preserving Operator Inference (EP-OpInf) approach, which imposes this structure on the learned reduced model via constrained optimization.
Numerical results using the viscous Burgers' and Kuramoto-Sivashinksy equation (KSE) demonstrate that EP-OpInf learns efficient and accurate reduced models that retain this energy-preserving structure.
\end{abstract}

\section{Nomenclature}

{\renewcommand\arraystretch{1.0}
\noindent\begin{longtable*}{@{}l @{\quad=\quad} p{5.5in}@{}}
$n$ & full model state dimension \\
$r$ & reduced model state dimension \\
$\mathbf A, \hat{\mathbf A}$ & full and reduced linear operators in $\mathbb{R}^{n\times n}$, $\mathbb{R}^{r\times r}$\\
$\mathbf H, \hat{\mathbf H}$ & full and reduced quadratic operators in $\mathbb{R}^{n\times n^2}$, $\mathbb{R}^{r\times r^2}$\\
$\mathbf F, \hat{\mathbf F}$ & full and reduced quadratic operators in $\mathbb{R}^{n\times n(n+1)/2}$, $\mathbb{R}^{r\times r(r+1)/2}$, containing unique columns of $\mathbf H, \hat{\mathbf H}$\\
$K$ & total number of data points \\
$t$ & time variable \\
$\mathbb R^{m\times n}$ & Euclidean space with dimension $m$-by-$n$ \\
$\mathbf V_r$ & POD basis in $\mathbb R^{n\times r}$ \\
$\mathbf x, \hat{\mathbf x}$ & full and reduced system states for linear-quadratic ODE in $\mathbb R^n,~\mathbb R^r$\\
$\mathbf X, \hat{\mathbf X}$ & full and reduced state snapshot data matrices in $\mathbb R^{n\times K},~\mathbb R^{r\times K}$ \\
$\hat {\mathbf X}_\otimes$ & reduced state quadratic data matrix in $\mathbb R^{r^2\times K}$ \\
$\dot{\mathbf X}, \dot{\hat{\mathbf{X}}}$  & full and reduced state time derivative data matrix in $\mathbb R^{n\times K}$, $\mathbb{R}^{r\times K}$\\
$\omega$ & spatial variable \\
$\otimes$ & Kronecker product \\
$\langle \cdot, \cdot \rangle$ & Euclidean inner product \\
$\|\cdot \|, \|\cdot\|_F$ & Euclidean and Frobenius norm \\
$\text{vech}(\cdot)$ & half-vectorization \\
$\text{Cov}[\cdot,\cdot]$ & covariance \\
$\text{Var}[\cdot]$ & variance
\end{longtable*}}

\section{Introduction}
\lettrine{M}{any-query} computations, which require repeated evaluations of a predictive model, arise in many different engineering contexts, including inverse problems, state estimation, uncertainty propagation, design, optimization, and control. In these many-query settings, high-dimensional models, such as those emerging from the spatial discretization of a PDE, can be prohibitively expensive to evaluate. Projection-based model reduction addresses this challenge by deriving an approximate model with a smaller dimension by projecting the original model governing equations onto a low-dimensional subspace  \cite{Benner_project_2015}, thereby reducing computational cost. Model reduction methods have enabled fast and cheap predictive simulations in many practical applications such as turbulent aerodynamic flow, weather forecasting, and large-scale structural analysis \cite{Rowley_romflow_2017,Horenko_weather_2008,Besselink_structure_2013}.

Model reduction methods are divided into intrusive and non-intrusive approaches. Intrusive methods are applicable
when there is direct access to both the model and its underlying computational code; this access is used to \textit{intrusively} compute reduced operators of the projected governing equations. Intrusive strategies include system theoretic methods like balanced truncation \cite{Moore_bt_1981} and rational interpolation methods \cite{grimme1997krylov,Antoulas2009-sj,Antoulas_interpolary_2010}, as well as sampling-based methods such as Proper Orthogonal Decomposition (POD) \cite{lumley1967structure,sirovich1987turbulence,berkooz1993proper}, Balanced Proper Orthogonal Decomposition \cite{Willcox_bpod_2002}, and reduced basis methods \cite{Peterson_rb_1989,Veroy_certified_2005,hesthaven2016certified}. Because intrusive reduced models are constructed by projecting the original governing equations onto a reduced subspace, intrusive methods can preserve some types of mathematical structure in the original model (stability, symmetries, invariances, etc.) that reflect physical principles governing the original dynamics. 

However, in many real-world systems, we encounter cases where we can compute state trajectories from given initial conditions but lack access to the full model operators or the code itself. This limitation necessitates the development of \textit{non-intrusive} model reduction methods which derive the reduced model operators from data rather than via intrusive projection. Non-intrusive methods which infer reduced matrix operators include Dynamic Mode Decomposition \cite{schmid_dmd_2010}, Operator Inference (OpInf) \cite{peherstorfer2016data,Qian_reduced_2022,mcquarrie_data-driven_2021}, and Lift \& Learn \cite{qian_lift_2020}. Several works have also developed approaches that use deep neural networks to learn low-dimensional dynamics~\cite{Swischuk_projection_2019,Bhattacharya_model_2021,Hesthaven_nn_2018,Gao_deep_2020}. An alternative family of methods uses sparse regression to identify from data the terms of a system's governing equations~\cite{Brunton_discover_2016,Rudy_pde_2017}. These approaches, leveraging available data and some understanding of system dynamics and structures, make surrogate modeling feasible for problems where the intrusive access to the model code is unavailable. However, because non-intrusive reduced model operators are not obtained from direct projection of the full model operators, they often fail to retain mathematical structure that the intrusive operators would inherit. This can yield reduced models that may lack accuracy, robustness, and generalizability \cite{baddoo_physics-informed_2023}. There is therefore a need to develop non-intrusive model reduction methods that preserve these types of structure.

Many structure-preserving reduced model learning methods have been developed, each incorporating different physical structures into their model learning process in different ways. Among these is a physics-informed dynamic mode decomposition method introduced in \cite{baddoo_physics-informed_2023}, which learns linear operators via a manifold optimization to preserve certain structures including circulant and symmetric structures. However, this method is applicable only to linear operators. For Hamiltonian systems, the works~\cite{sharma_hamiltonian_2022,Gruber_nchopinf_2023} constrain the Operator Inference approach to retain Hamiltonian structure.  
Further contributing to this field, \cite{sharma_lagrange_2022} developed a structure-preserving OpInf method that ensures the inferred system operators are Lagrangian and can preserve symmetry or symmetric positive definite structure. These structure-preserving non-intrusive methods have been demonstrated on applications ranging from the linear wave equation to Schr\"odinger's equation to soft robots \cite{baddoo_physics-informed_2023,sharma_hamiltonian_2022,Gruber_nchopinf_2023,sharma_lagrange_2022}.

In this paper, we consider quadratic systems in which the quadratic term of the governing equations is energy-preserving: this means that the contribution of the quadratic term to the time derivative of the system energy is zero~\cite{schlegel_noack_2015}. Several works have already incorporated this energy-preserving structure into different non-intrusive model learning methods. For instance, \cite{kaptanoglu2021promoting} modified the sparse regression approach \cite{Brunton_discover_2016,Rudy_pde_2017} to promote learning quadratic terms with this structure \cite{Kaptanoglu_magneto_2021}. The work \cite{goyal_gurantee_2023} proposes a method for inferring reduced operators that leads to local and global stability guarantees; the energy-preserving structure of the quadratic term is enforced through a specific parametrization of the reduced quadratic operator. In contrast, the present work focuses solely on enforcing the energy-preserving structure of the quadratic term through a constrained optimization; providing stability guarantees for models learned in this way is a direction for future research.



We introduce Energy-Preserving OpInf (EP-OpInf), an OpInf method that enforces energy-preserving structure in the learned quadratic operator via a constrained optimization. Numerical experiments are conducted for the viscous Burgers' equation and the Kuramoto-Sivashinsky equation (KSE). Numerical results demonstrate that EP-OpInf learned models achieve accuracy similar to that of the unconstrained data-fit model learned via standard OpInf while non-intrusively retaining the energy-preserving structure that intrusive models inherit. This lays the foundation for applying stability analyses developed for intrusive reduced models to the EP-OpInf learned models. 

The remainder of the paper is organized as follows. Section \ref{sec:2-background} introduces the quadratic systems that we consider and provides an overview of both intrusive model reduction via POD and the standard non-intrusive OpInf method. Section \ref{sec:3-eng-preserve-opinf} then presents our proposed EP-OpInf method and discusses some implementation considerations. Section \ref{sec:4-experiment} demonstrates the efficacy of our proposed method through a series of numerical experiments. Conclusions are presented in Section \ref{sec:5-conclusion}.

\section{Background}\label{sec:2-background}

In this section, we introduce the quadratic system with energy-preserving quadratic term that we consider (Section \ref{subsec:2-ep-background}), the POD method for model reduction (Section \ref{subsec:2-pod}), and the standard OpInf method for learning reduced models non-intrusively (Section \ref{subsec:2-opinf}).

\subsection{Energy Preserving Quadratic Nonlinearities in PDEs}\label{subsec:2-ep-background}
We consider an $n$-dimensional ordinary differential equation (ODE) which is linear and quadratic in the state $\mathbf x$. Such an ODE often arises from spatially discretizing a PDE and is given by
\begin{equation}\label{eqn:quadratic-system}
    \dot{\mathbf x}(t) = \mathbf A \mathbf x(t) + \mathbf H \left(\mathbf x(t) \otimes \mathbf x(t)\right) 
\end{equation}
where $\mathbf x(t) \in \mathbb{R}^n$ is the system state vector over $t \in [0, T_{\text{final}}]$, and $\otimes$ denotes the Kronecker product. The operators $\mathbf A \in \mathbb{R}^{n\times n}$ and $\mathbf H \in \mathbb{R}^{n\times n^2}$ are the linear and quadratic operators, respectively. In our setting, $n$ is large, so simulating the system~\eqref{eqn:quadratic-system} is computationally expensive.

The quadratic operator $\mathbf H$ is called `energy-preserving' if
\begin{equation}\label{eqn:ep-property}
    \langle \mathbf x, \mathbf H (\mathbf x \otimes \mathbf x)\rangle = \mathbf x^\top \mathbf H(\mathbf x \otimes \mathbf x) = 0, \qquad \text{ for all } \mathbf x \in \mathbb R^n.
\end{equation}
This condition is derived by setting the quadratic term in the time derivative of the energy, $\frac12\|\mathbf{x}\|^2$, to zero. Note that this condition can be generalized to energies defined by weighted norms or inner products; our EP-OpInf approach also would generalize to such settings following the weighted OpInf approach of \cite{Qian_reduced_2022}. Such energy-preserving quadratic terms commonly arise from nonlinear convection terms in the governing equations for fluid problems, including the Navier-Stokes equations. Some types of stability guarantees are possible for models with energy-preserving quadratic nonlinearities \cite{schlegel_noack_2015}, making it desirable to enforce this property in learning reduced models that will be used for design and control.

\subsection{Proper Orthogonal Decomposition}\label{subsec:2-pod}
Proper orthogonal decomposition \cite{lumley1967structure, sirovich1987turbulence, berkooz1993proper} originated from the analysis of turbulent flows in aerodynamics, and it has become one of the most widespread projection-based model reduction methods. POD reduces the model by projecting it onto a reduced subspace defined to be the span of basis vectors that optimally represent a set of simulation or experimental data.

In POD, we begin by collecting snapshots of state trajectory time series data by simulating the original full model ODE \eqref{eqn:quadratic-system} with $K$ timesteps. We define the state snapshot data matrix as follows: 
\begin{equation}\label{eqn:state-input-snapshot-matrix}
    \boldsymbol{\mathbf X} = \begin{bmatrix} 
        | & | &  & | \\
        \boldsymbol{\mathbf x}(t_1) & \boldsymbol{ \mathbf x}(t_2) & \cdots & \boldsymbol{\mathbf x}(t_K) \\
        | & | &  & | \\
    \end{bmatrix} \in \mathbb{R}^{n\times K} ~.
\end{equation}
More generally, the state snapshot matrix can contain state data from multiple simulations, e.g., from different initial conditions or using different parameters.
Let $\mathbf{X} = \mathbf{V\Sigma W}^\top$ denote the singular value decomposition of the state snapshot matrix~\eqref{eqn:state-input-snapshot-matrix}. To reduce the dimension of the large-scale model, we denote by $\mathbf{V}_r\in \mathbb R^{n\times r}$ the first $r \ll n$ columns of $\mathbf V$; this is called the \textit{POD basis}. Then, we approximate the state $\mathbf{x}$ in the subspace spanned by the POD basis, $\mathbf x \approx \mathbf V_r \hat{\mathbf x}$ where $\hat{\mathbf x}\in\mathbb{R}^r$ is called the \textit{reduced state}. If we substitute this approximation into \eqref{eqn:quadratic-system} and enforce the Galerkin orthogonality condition that the approximation residual be orthogonal to the span of $\mathbf V_r$, we arrive at a POD-Galerkin reduced model of the form
\begin{equation}\label{eqn:reduced-quadratic}
    \dot{\hat{\mathbf x}}(t) = \hat{\mathbf A}\hat{\mathbf x}(t) + \hat{\mathbf{H}}(\hat{\mathbf{x}}(t) \otimes \hat{\mathbf{x}}(t)),
\end{equation}
where the reduced operators are $\hat{\mathbf{A}} = \mathbf{V}^\top_r \mathbf{AV}_r \in \mathbb{R}^{r\times r}$ and  $\hat{\mathbf{H}} = \mathbf{V}^\top_r \mathbf{H}(\mathbf{V}_r \otimes \mathbf{V}_r) \in \mathbb R^{r\times r^2}$.
Since the size of the reduced operators scales with the reduced dimension $r$ and not the original dimension $n$, we are able to evaluate the reduced model efficiently in many-query computations arising in design and control. The computation of the reduced operators in traditional POD-Galerkin model reduction is \textit{intrusive}, that is, we require access to the full operators $\mathbf{A}$ and $\mathbf{H}$ in order to compute their projections.  The energy-preserving property of the original matrix $\mathbf{H}$ \eqref{eqn:ep-property} is inherited in the reduced matrix $\hat{\mathbf{H}}$ through this intrusive projection process.

However, cases exist where access to the original operators is limited or impossible, but we know the model structure, such as linear, quadratic, etc., based on the application. In such situations, Operator Inference~\cite{peherstorfer2016data}, described in Section \ref{subsec:2-opinf}, can \textit{non-intrusively} learn reduced operators solely from state trajectory and input data. 

\subsection{Operator Inference}\label{subsec:2-opinf}
The goal of Operator Inference is to non-intrusively obtain a reduced model of the form~\eqref{eqn:reduced-quadratic}. To do so, we will fit reduced operators $\hat{\mathbf{A}}$ and $\hat{\mathbf{H}}$ to the reduced data in a least-squares sense. In addition to the state trajectory data \eqref{eqn:state-input-snapshot-matrix}, we also require paired state time derivative data: 
\begin{align*}
    \dot{\mathbf{X}} = \begin{bmatrix}
        | & | & & | \\
        \dot{\mathbf{x}}(t_1) & \dot{\mathbf{x}}(t_2) & \cdots & \dot{\mathbf{x}}(t_K) \\
        | & | & & | 
    \end{bmatrix}\in\mathbb{R}^{n\times K} \quad \text{ where ~~} 
    \dot{\mathbf x}(t_i)=\mathbf A\mathbf x(t_i)+ \mathbf H (\mathbf x(t_i) \otimes \mathbf x(t_i))~.
\end{align*}
This time derivative data can come directly from the simulation of the high-dimensional ODE or can be approximated numerically from the snapshot data. We use the POD basis $\mathbf{V}_r$ to compute reduced state and time derivative data as follows: let $\hat{\mathbf{x}}_i = \mathbf{V}_r^\top\mathbf{x}(t_i)$, and $\dot{\hat{\mathbf{x}}}_i = \mathbf{V}_r^\top\dot{\mathbf{x}}(t_i)$ for $i=1,\ldots,K$. Then, define
\begin{equation}
    \hat{\textbf X} = \begin{bmatrix}
        | & | &  & | \\
        \boldsymbol{\hat{\mathbf x}}_1 & \boldsymbol{ \hat{\mathbf x}}_2 & \cdots & \boldsymbol{\hat{\mathbf x}}_K \\
        | & | &  & | \\
    \end{bmatrix} \in \mathbb R^{r\times K}, \qquad \text{and} \qquad 
    \dot{\hat{\textbf X}} = \begin{bmatrix}
        | & | &  & | \\
        \dot{\hat{\mathbf x}}_1 & \dot{\hat{\mathbf x}}_2 & \cdots & \dot{\hat{\mathbf x}}_K \\
        | & | &  & | \\
    \end{bmatrix} \in \mathbb R^{r\times K}.
\end{equation}
Additionally, we define the matrix formed by the quadratic terms of the state data
\begin{equation}\label{eqn:H-square-matrix}
    \hat{\textbf X}_\otimes = \begin{bmatrix}
        | & | &  & | \\
        (\hat{\mathbf x}_1 \otimes \hat{\mathbf x}_1) & (\hat{\mathbf x}_2 \otimes \hat{\mathbf x}_2) & \cdots & (\hat{\mathbf x}_K \otimes \hat{\mathbf x}_K) \\
        | & | &  & | \\
    \end{bmatrix} \in \mathbb R^{r^2\times K}.
\end{equation}
This allows us to formulate the following minimization for finding the reduced operators $\hat{\mathbf{A}}$ and $\hat{\mathbf{H}}$:   
\begin{equation}\label{eqn:optimization1}
    \textbf{Standard OpInf}: \qquad 
    \min_{\hat{\mathbf{A}}\in\mathbb R^{r\times r},~\hat{\mathbf{H}}\in\mathbb R^{r\times r^2}} \sum_{i=1}^K \left\| \dot{\hat{\mathbf x}}_i - \hat{\mathbf{A}}\hat{\mathbf{x}}_i - \hat{\mathbf H}(\hat{\mathbf x}_i \otimes \hat{\mathbf x}_i) \right \|^2_2 = \min_{\mathbf{O}\in\mathbb R^{r\times(r+r^2)}} \|\mathbf{D}\mathbf{O}^\top - \dot{\hat{\mathbf X}}\|_F^2,
\end{equation}
where $\mathbf D = [\hat{\mathbf X}^\top, ~\hat{\mathbf X}_\otimes^\top] \in \mathbb R^{K\times(r+r^2)}$ and $\mathbf O = [\hat{\mathbf A}, ~\hat{\mathbf H}] \in \mathbb R^{r\times(r+r^2)}$. We note that the redundancy in the Kronecker product means that the minimizer $\hat{\mathbf{H}}$ of \eqref{eqn:optimization1} is non-unique; this non-uniqueness can be eliminated by requiring its entries to satisfy a symmetry condition (see Section~\ref{subsec:3-implementation-considerations} for a discussion on implementation). The optimization \eqref{eqn:optimization1} is a linear least-squares problem and can be solved efficiently using standard numerical linear algebra packages. In fact, the problem~\eqref{eqn:optimization1} decouples into $r$ independent least-squares problems, one for each row of the reduced operators with $r+r^2$ degrees of freedom each.

\section{Energy-Preserving Operator Inference}\label{sec:3-eng-preserve-opinf}
Section \ref{subsec:3-energy-preserving-constraint} introduces our Energy-Preserving Operator Inference method. Section \ref{subsec:3-implementation-considerations} discusses efficient implementation of the method.

\subsection{Energy-Preserving Constraint}\label{subsec:3-energy-preserving-constraint}
The Operator Inference approach in the preceding section finds operators that best fit the reduced data. However, it does not necessarily find a quadratic operator $\mathbf{H}$ with the energy-preserving property~\eqref{eqn:ep-property}. Thus, we now introduce Energy-Preserving Operator Inference (EP-OpInf), which learns a quadratic reduced model with an energy-preserving quadratic operator. 

We first re-write the energy-preservation property~\eqref{eqn:ep-property} as a constraint on the entries of $\hat{\mathbf{H}}$. Let $\hat{\mathbf H}[l,m]$ correspond to an element of $\hat{\mathbf H}$ at the row index $l$ and column index $m$, and define 
\begin{equation}
    \hat h_{ijk} = \hat{\mathbf H}[i,r(k-1)+j].
\end{equation}
We assume without loss of generality that the entries of $\hat{\mathbf H}$ are symmetric in the last two indices, i.e., $\hat h_{ijk} = \hat h_{ikj}$ (in Section \ref{subsec:3-implementation-considerations} we discuss how we enforce this in implementation). Then, the energy-preserving condition can be rewritten in terms of the entries of $\hat{\mathbf H}$ as \cite{kaptanoglu2021promoting,schlegel_noack_2015}:
\begin{equation}\label{eqn:energy-preservation-H}
    \hat h_{ijk} + \hat h_{jik} + \hat h_{kij} = 0, \quad  1 \leq i,j,k \leq r ~ .
\end{equation}
Due to the symmetry of the last two indices of $\hat h_{ijk}$, this defines a total of $\frac{r(r+1)(r+2)}{6}$ unique constraints. We emphasize that submatrices of size $r'<r$ extracted from the $\hat{\mathbf H}$ corresponding to $i,j,k=1,\ldots,r'$ also satisfy the energy-preserving condition, since the constraint \ref{eqn:energy-preservation-H} is true for all index combinations. 

To impose this energy-preserving structure on the operator, we propose EP-OpInf. For this method, we incorporate the constraint~\eqref{eqn:energy-preservation-H} into optimization~\eqref{eqn:optimization1} and formulate a constrained minimization as follows:
\begin{equation}\label{eqn:ephec-opinf}
    \textbf{EP-OpInf}: \qquad 
    \min_{\mathbf{O}\in\mathbb R^{r\times(r+r^2)}} \|\mathbf{D}\mathbf{O}^\top - \dot{\hat{\mathbf{X}}}\|_F^2 \quad \text{ subject to } \quad \hat h_{ijk} + \hat h_{jik} + \hat h_{kij} = 0, \quad  1 \leq i,j,k \leq r ~ .
\end{equation}
We note that the original Operator Inference approach \eqref{eqn:optimization1} was an unconstrained linear least-squares problem that could be efficiently decoupled into individual rows of the operators and solved with standard numerical linear algebra packages. 
With the addition of our constraints, the minimization can no longer be decoupled into rows, and the minimization must be solved using constrained optimization solvers. There is thus a trade-off between the simplicity of implementing the Operator Inference approach and preserving the energy-preserving structure of the quadratic term.

\subsection{Implementation Considerations}\label{subsec:3-implementation-considerations}
For the efficient computation of the optimization \eqref{eqn:ephec-opinf}, we reformulate the quadratic operator and associated data matrix constructed by the Kronecker product. The code implementations for this reformulation are available online\footnote{GitHub repository: \href{https://github.com/smallpondtom/LiftAndLearn.jl}{https://github.com/smallpondtom/LiftAndLearn.jl}}. We note that the matrix $\hat{\mathbf{X}}_\otimes$ \eqref{eqn:H-square-matrix} contains redundant rows due to the symmetry of the Kronecker product. Thus, the solution of the minimization~\eqref{eqn:ephec-opinf} is not unique unless we enforce symmetry in the last two indices of the learned operator $\hat{\mathbf{H}}$. To do so, we will re-express the optimization in terms of an operator $\hat{\mathbf{F}}\in\mathbb{R}^{r\times r(r+1)/2}$, which removes the redundant columns of $\hat{\mathbf{H}}$. This has the added advantage of reducing the dimension of the optimization. This operator $\hat{\mathbf{F}}$ is defined as follows. Let $\hat f_{ijk}$ denote the entry of the $i$-th row in column $(r-k/2)(k-1)+j$, that is:
\begin{equation}
    \hat{f}_{ijk} = \hat{\mathbf F}[i,(r-k/2)(k-1)+j] ~ , \quad \text{for }~~i = 1,\ldots, r; ~j = 1,\ldots, r; ~k = 1,\ldots, j.
\end{equation}
Then, the entries of $\hat{\mathbf{H}}$ and $\hat{\mathbf{F}}$ can be related as follows:
\begin{equation}
    \hat{f}_{ijk} = \begin{cases}
        \hat{h}_{ijk}  & \textrm{ if } j = k \\
        \hat{h}_{ijk} + \hat{h}_{ikj}  & \textrm{ otherwise }
    \end{cases}, \quad i = 1,\ldots, r; ~j = 1,\ldots, r; ~k = 1,\ldots, j ~ .
\end{equation}
With this relation and considering the permutations of indices, we rewrite the quadratic energy-preserving constraint as
\begin{equation}
    \delta_{jk}\hat{f}_{ijk} + \delta_{ik}\hat{f}_{jik} + \delta_{ij}\hat{f}_{kij} = 0 ~ , \quad \textrm{where} \quad \delta_{jk} = \begin{cases}
        1 & \textrm{ if }~ j = k\\
        1/2 & \textrm{ otherwise }
    \end{cases}, \quad i = 1,\ldots, r; ~j = 1,\ldots, i; ~k = 1,\ldots, j ~ .
\end{equation}
Denote $\mathbf x^{[2]} = \text{vech}(\mathbf{xx}^\top)$ (`half-vectorization' of $\mathbf{xx}^\top$), to be obtained by eliminating the supradiagonal elements of $\mathbf{xx}^\top$ and vectorizing it, i.e., $x^{[2]}_{(r-k/2)(k-1)+j} = x_j x_k$ \cite{magnus1980Elimination}. Using this operation, let
\begin{equation}
    \hat{\textbf X}^{[2]} = \begin{bmatrix}
        | & | &  & | \\
        \hat{\mathbf x}_1^{[2]} & \hat{\mathbf x}_2^{[2]} & \cdots & \hat{\mathbf x}_K^{[2]} \\
        | & | &  & | \\
    \end{bmatrix} \in \mathbb R^{(r(r+1)/2)\times K}.
\end{equation}
Then, we can reformulate the optimization \eqref{eqn:ephec-opinf} with matrices $\mathbf{\widetilde D} = [\hat{\mathbf X}^\top, \hat{\mathbf X}^{[2]\top}] \in \mathbb R^{r\times (r+r(r+1)/2)}$ and $\mathbf{\widetilde O} = [\hat{\mathbf A}, \hat{\mathbf F}] \in \mathbb R^{r\times(r+r(r+1)/2)}$ as 
\begin{equation}
    \min_{{\mathbf{\widetilde O}} \in \mathbb R^{r\times(r+r(r+1)/2)}} \| \mathbf{\widetilde D}\mathbf{\widetilde O^\top}- \mathbf R \|^2_F \quad \text{ subject to } \quad \delta_{jk}\hat{f}_{ijk} + \delta_{ik}\hat{f}_{jik} + \delta_{ij}\hat{f}_{kij} = 0 ~ , \quad i = 1,\ldots, r; ~j = 1,\ldots, i; ~k = 1,\ldots, j ~ .
\end{equation}
This reformulation decreases the degrees of freedom in the least squares problem associated with the quadratic operator from $r^2$ to $r(r+1)/2$, and also reduces the number of constraints in the minimization problem, making it more efficient than a na\"ive implementation of~\eqref{eqn:ephec-opinf}.

\section{Numerical Experiments}\label{sec:4-experiment}
This section demonstrates the proposed EP-OpInf method through its application to the viscous Burgers' equation (Section \ref{subsec:4-burgers}) and the Kuramoto-Sivashinksy equation (Section \ref{subsec:4-kse}).

\subsection{Viscous Burgers' Equation}\label{subsec:4-burgers}
In our first numerical experiment, we consider the viscous Burgers' equation, a simplified fluid dynamics model that describes the motion of a one-dimensional viscous fluid \cite{burgers1948mathematical}. It is given as:
\begin{equation}\label{eqn:burgers}
    \frac{\partial x}{\partial t}(\omega,t) = \mu \frac{\partial^2 x}{\partial \omega^2}(\omega, t) - x(\omega,t)\frac{\partial x}{\partial \omega}(\omega,t)~,
\end{equation}
where $x(\omega,t)$ is the velocity of the fluid, $t$ is time, $\omega$ is the spatial coordinate, and $\mu$ is the kinematic viscosity coefficient.  The viscous Burgers' equation can be derived from the Navier-Stokes equations under certain assumptions, including the assumption of one-dimensional flow and the neglect of pressure gradient driving the flow. The Burgers' equation is suitable for testing our EP-OpInf method since it is expressed in the form of \eqref{eqn:quadratic-system} where the energy-preserving quadratic nonlinearity holds \cite{aref_note_1984}.
\begin{figure}[h!]
    \centering
    \includegraphics[width=0.49\textwidth]{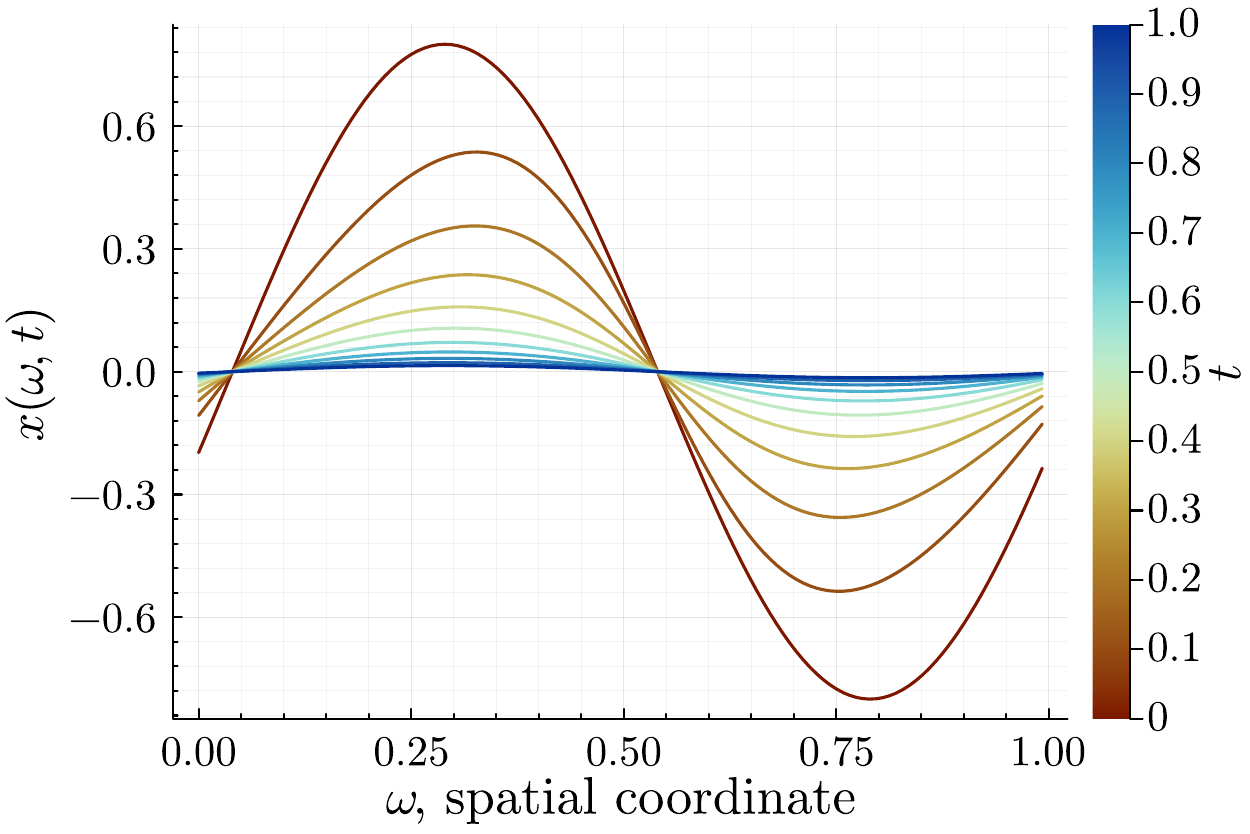}\hfill
    \includegraphics[width=0.49\textwidth]{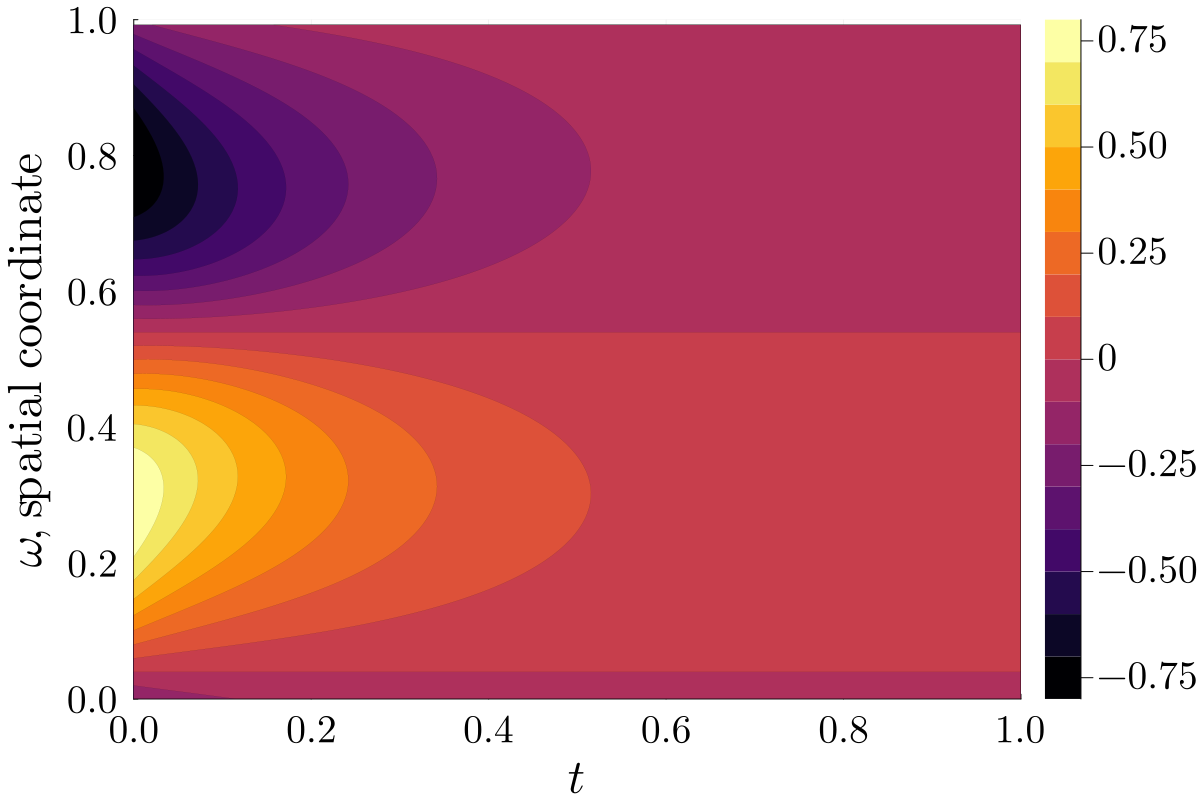}
    \caption{Visualization of Burgers' equation state evolution and flow field with $\mu = 0.1$ and reduced dimension of 10 with initial condition $x(\omega,0) = 0.8\sin(2\pi\omega -0.125)$.}
    \label{fig:burger-visual}
\end{figure}

\subsubsection{Setup}\label{subsubsec:4-burgers-setup}
We solve Burgers' equation on the unit domain with periodic boundaries, $\omega \in [0, 1)$, for $t \in [0, T]$, and for $\mu=0.1$. We use a sinusoidal initial condition $x(\omega,0) = A\sin(2\pi f  \omega + \phi)$ with three varied parameters: the amplitude $A$, frequency $f$, and phase $\phi$. We collect training data at $A=[0.8, 0.9, 1.0, 1.1, 1.2]$, $f=[1,2,3]$, and $\phi=[-0.25, -0.125, 0, 0.125, 0.25]$. In total, $5\times 3\times 5=75$ different parameterized initial conditions define the training data set.

This system is discretized on a uniform spatial grid with a grid spacing $\Delta\omega=2^{-7}$, leading to an ODE of the form \eqref{eqn:quadratic-system}. With this grid size, the state dimension is $n = 128$. The ODE is solved using a semi-implicit Euler scheme with a time step size of $\Delta t=10^{-4}$ and a final time of $T = 1$. Figure~\ref{fig:burger-visual} shows the behavior of the state for sample parameter values. 

In the model learning phase, data at every 100\textsuperscript{th} timestep are collected from the simulation described above, and these data are used to construct the data snapshot matrices and train an inferred reduced model.
To assess the trained operators, we conduct two tests where we generate additional data with different initial conditions:
\begin{itemize}
    \item \textbf{Test 1} (\textbf{Interpolation}): 50 different initial conditions generated uniformly within the training region: \\
    $A\in [0.8,  1.2]$, $f\in \{1,2,3\}$, and $\phi\in[-0.25,  0.25]$
    \item \textbf{Test 2} (\textbf{Extrapolation}): 50 different initial conditions generated uniformly outside the training region: \\
    $A \in [0.5,0.8) \cup (1.2,1.5]$, $f \in \{4,5,6\}$, and $\phi \in [-0.5,-0.25) \cup (0.25, 0.5]$.
\end{itemize}\vspace{0.3cm}

\begin{figure}[t!]
    \centering
    \includegraphics[width=0.47\textwidth]{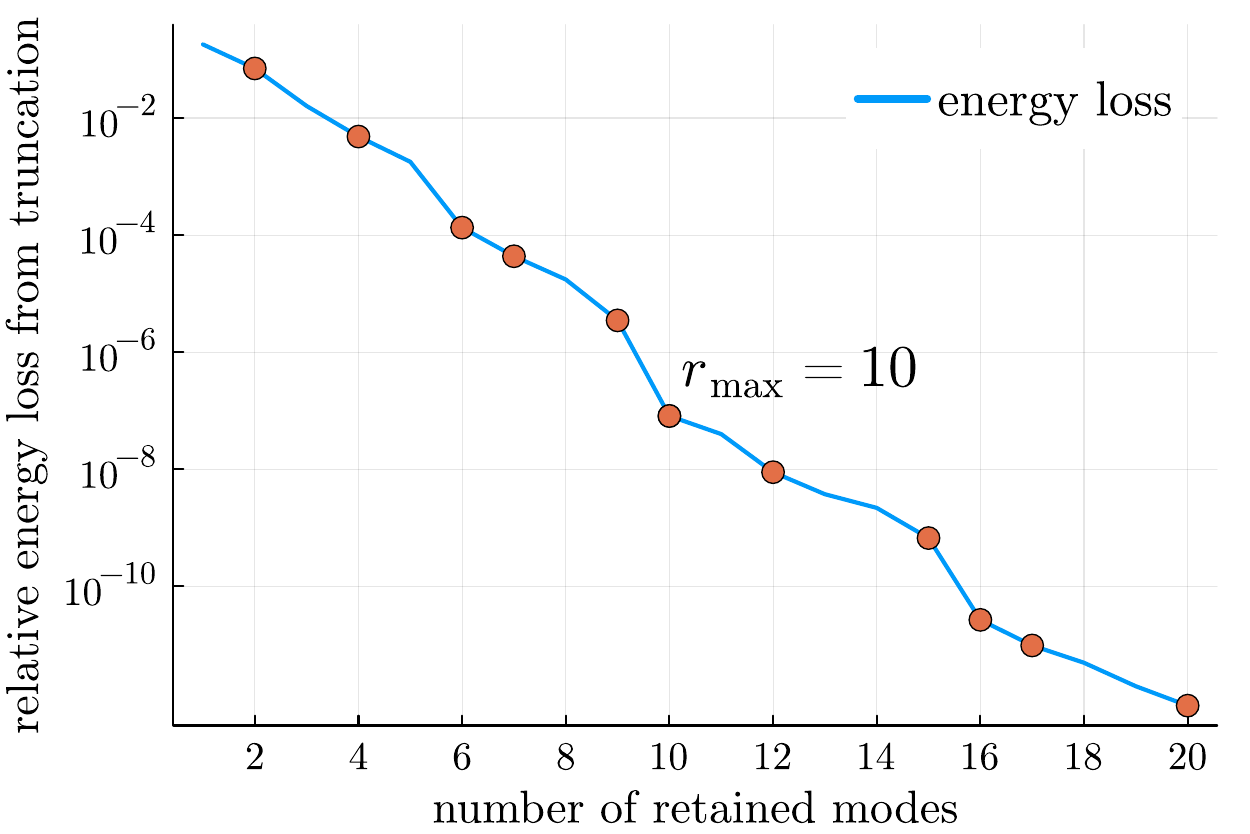}
    \caption{Relative energy lost from truncated modes for viscous Burgers' equation training data.}
    \label{fig:burgers-eng-spectrum}
\end{figure}

To determine the number of POD modes to retain, we compute the relative energy lost by a POD basis of size $r$, defined to be
\begin{equation}\label{eqn:energy-spectrum}
    1 - \frac{\sum_{i=r+1}^n \sigma_i^2}{\sum_{i=1}^n \sigma_i^2}~.
\end{equation}
\noindent In this equation, $\sigma_i$ represents the singular values of the snapshot data matrix. The energy spectrum for the Burgers' equation in Figure \ref{fig:burgers-eng-spectrum} shows that when $r_{\text{max}}=10$ we capture 99.99999\% of the energy in the data. We therefore set $r_{\text{max}} = 10$ in our numerical experiments.
For both standard OpInf and EP-OpInf, the reduced operators for dimensions $r < r_{\text{max}}$ are retrieved by extracting the appropriate submatrices from $r_{\text{max}}$-dimensional model operators. As mentioned in Section \ref{subsec:3-energy-preserving-constraint}, the energy-preserving structure holds even for extracted operators.
 Once we train the inferred models, we use the learned models to reconstruct the training data by simulating them from the same initial condition that generated the training data. We report the relative state error of the reduced model reconstruction, defined as
\begin{equation}\label{eqn:rel-state-err}
    \text{average relative state error} = \frac{1}{J}\sum_{j=1}^J \frac{\left\|\mathbf X_j - \overline{\mathbf X}_j \mathbf V_r^\top\right\|^2_F}{\left \| \mathbf X_j \right\|_F^2}~,
\end{equation}
where $J$ is the total number of initial conditions and $\overline{\mathbf X}$ contains the trajectory data obtained by integrating the surrogate model defined by either the intrusive \eqref{eqn:reduced-quadratic} or inferred operators with either standard OpInf~\eqref{eqn:optimization1} or EP-OpInf \eqref{eqn:ephec-opinf}. 

Additionally, because a key goal of this work is to learn models whose quadratic terms have the energy-preserving property, we also report the following measure of the violation of the energy-preserving property:
\begin{equation}\label{eqn:constraint-violation}
    \text{energy-preserving constraint violation} = \sum_{1 \leq i , j , k \leq r} |\hat{h}_{ijk} + \hat{h}_{jik} + \hat{h}_{kji}|~,
\end{equation}
which is defined to be the sum of constraint residuals. 

\subsubsection{Results}\label{subsubsec:4-burgers-results}
\begin{figure}[h!]
    \centering
    \includegraphics[width=0.47\textwidth]{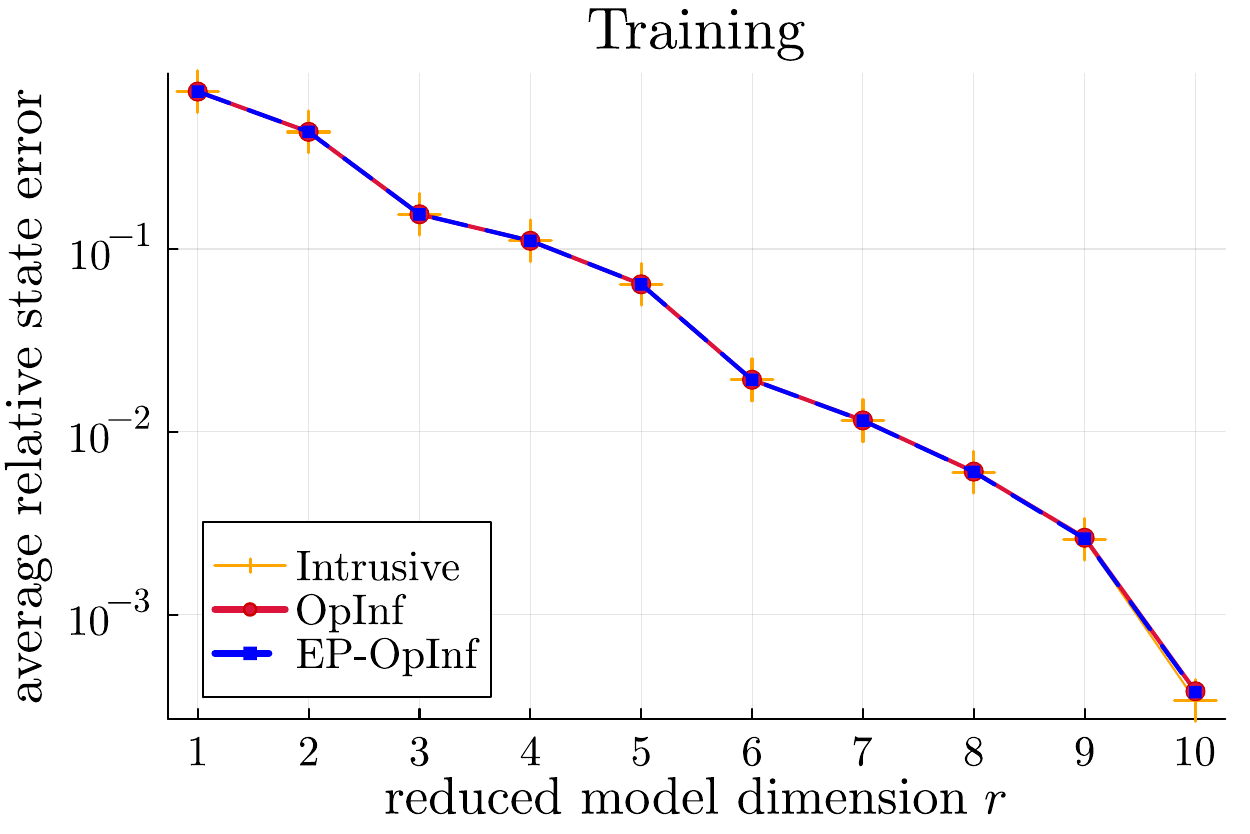}\hfill
    \includegraphics[width=0.47\textwidth]{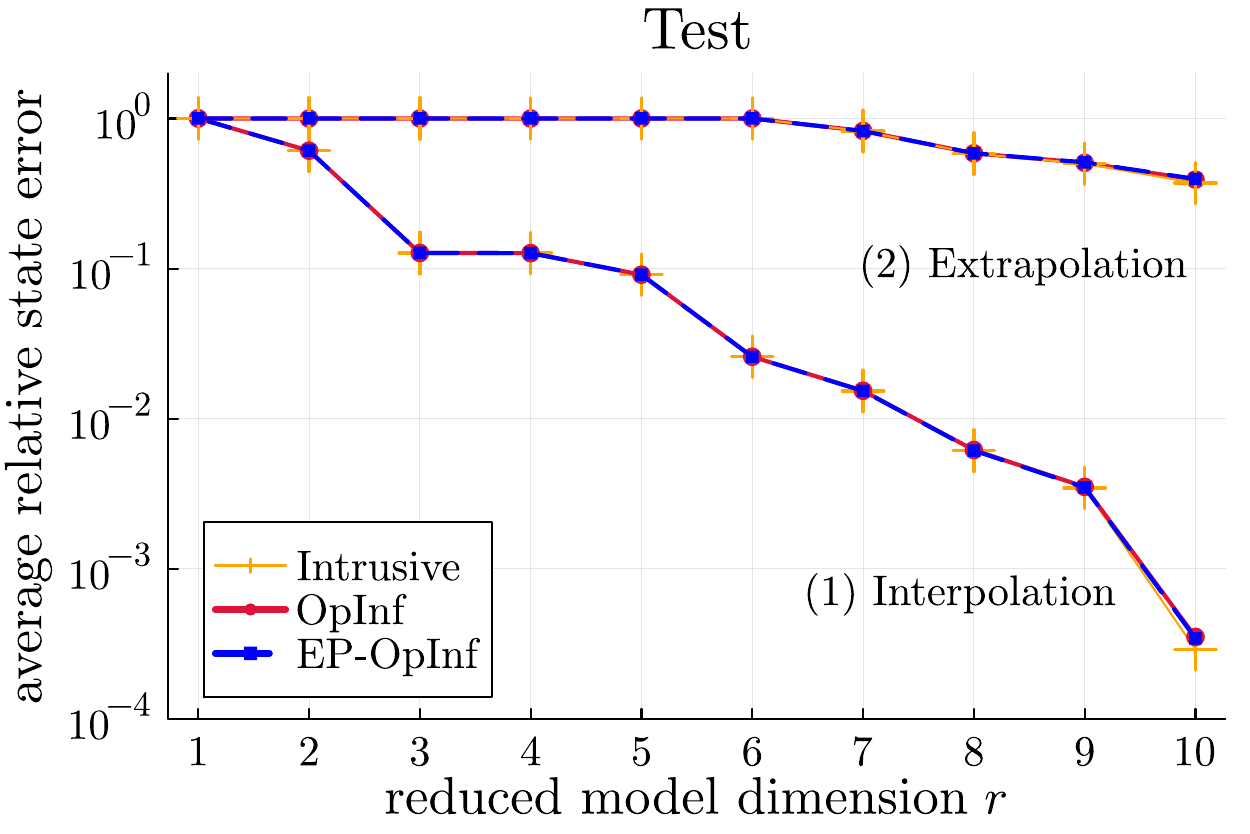}
    \caption{Comparison of intrusive, standard OpInf, and EP-OpInf methods of relative state error for training data \ref{eqn:rel-state-err} (left) and test cases (right).}
    \label{fig:burgers-rse}
\end{figure}

\begin{figure}[h!]
    \centering
    \includegraphics[width=0.47\textwidth]{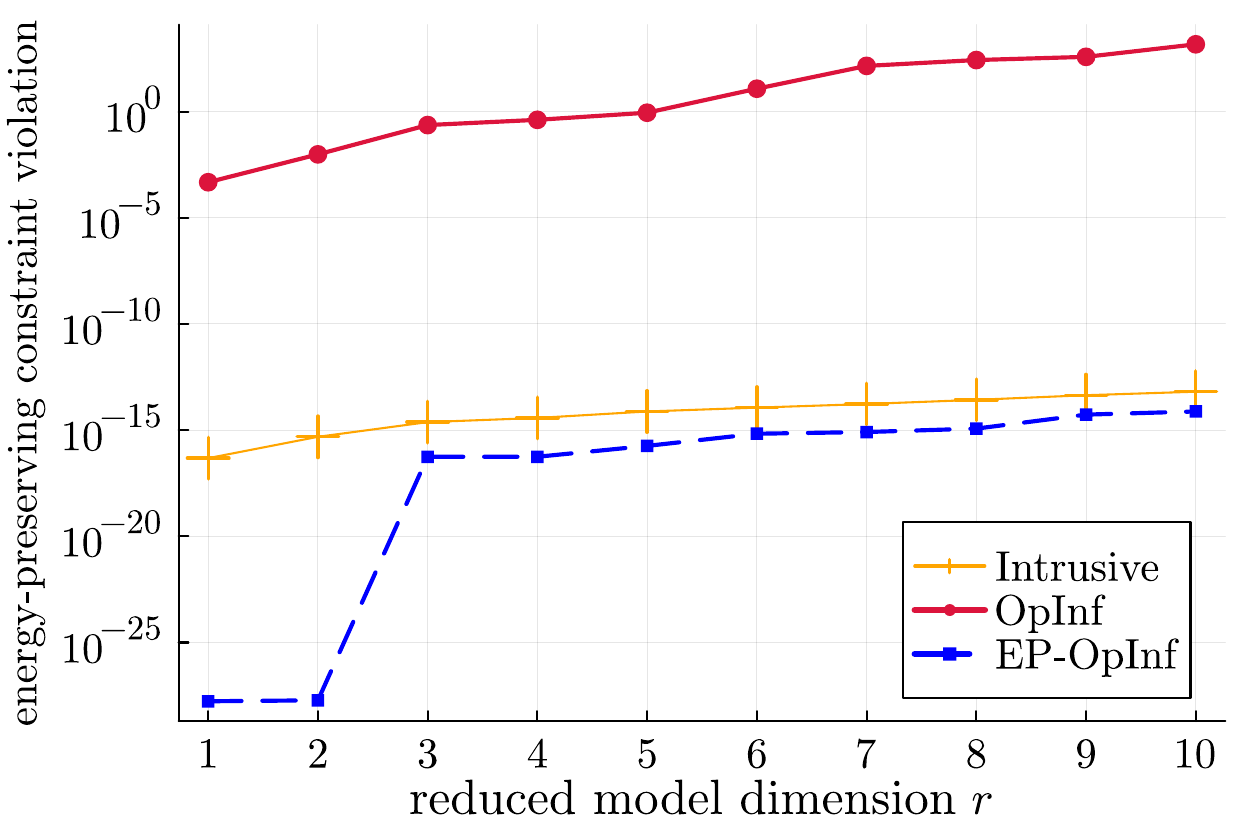}
    \caption{Constraint violation of all three methods for Burgers' equation over reduced dimensions $r$.}
    \label{fig:burgers-cv}
\end{figure}

The left plot in Figure \ref{fig:burgers-rse} shows the relative state error for the training data computed from the intrusive, standard OpInf, and EP-OpInf models. Figure \ref{fig:burgers-cv} plots the energy-preserving constraint violations vs.\ the size of the reduced model. From Figure \ref{fig:burgers-rse}, we see the accuracy of all methods in reconstructing the training trajectories is similar. From the constraint violations in Figure \ref{fig:burgers-cv}, it is evident that the intrusive and EP-OpInf models satisfy the energy-preserving condition up to machine precision, indicating that the quadratic operators of the two models possess the energy-preserving property. We emphasize that EP-OpInf does this \textit{non-intrusively}, whereas the intrusive method requires access to the original model code. Additionally, despite the learning process for EP-OpInf being constrained by the energy-preserving quadratic structure, EP-OpInf achieves a level of model accuracy that is very close to both that of the unconstrained standard OpInf method and the intrusive method.

The right plot in Figure \ref{fig:burgers-rse} plots the test errors for the intrusive, standard OpInf, and EP-OpInf models, for both the test data sets considered. For different initial conditions generated within the training region, the accuracy is comparable to that of the training data results shown in the plot on the left, with all three methods yielding similar levels of accuracy. Even if we extrapolate from the training region as in Test 2, the standard OpInf and EP-OpInf errors are on par with the intrusive method. We emphasize that both the standard and EP-OpInf models \textit{non-intrusively} achieve accuracy comparable to that of the intrusive model, and that EP-OpInf does this while enforcing the energy-preserving structure. 


\subsection{Kuramoto-Sivashinsky Equation}\label{subsec:4-kse}
Originally introduced in the context of flame propagation, the Kuramoto-Sivashinsky equation has since found applications in various fields, including fluid dynamics, pattern formation, and combustion theory \cite{kuramoto1978diffusion,sivashinksy1977nonlinear}. The equation is given by:
\begin{equation}\label{eqn:kse}
    \frac{\partial x}{\partial t}(\omega,t) + \frac{\partial^2 x}{\partial \omega^2}(\omega, t) + \mu \frac{\partial^4 x}{\partial \omega^4}(\omega, t) + x(\omega,t)\frac{\partial x}{\partial \omega}(\omega,t) = 0~.
\end{equation}
Here, $x(\omega,t)$ represents the fluid velocity, $t$ denotes time, $\omega$ is the spatial coordinate, and $\mu$ is the kinematic viscosity coefficient. Like the viscous Burgers' equation, the KSE satisfies the energy-preserving quadratic nonlinearity \cite{aref_note_1984}. However, the KSE poses challenges for model reduction. This is because the KSE models localized turbulence in which coherent spatial structures and temporal chaos coexist \cite{hyman_order_1986}.
The temporal chaos indicates the exponential separations of nearby state trajectories of the KSE solutions; thereby, even small errors in the initial condition can be amplified, making the long-time dynamics unpredictable \cite{Strogatz_chaos_2018}. This unpredictability increases the error when reconstructing the state trajectories from the reduced models.
In contrast, the coherent spatial structures indicate inertial manifolds, i.e., invariant manifolds containing global attractors for dissipative dynamical systems \cite{Constantin_manifolds_1989,Lu_databased_2017,cvitanovic_state_2010}. The solutions of KSE are attracted asymptotically to this manifold, and it is important for OpInf to capture this finite-dimensional manifold accurately. 

Several studies in the literature have explored model reduction approaches for the KSE. For instance, one study employed an intrusive Galerkin projection method with regularization to derive a reduced model and compared the state trajectories to the direct numerical simulations \cite{sabetghadam_regularize_2012}. Other studies implemented non-intrusive methods such as OpInf and NARMAX (Nonlinear AutoRegressive Moving Average with eXogenous input) \cite{Chen_narmax_1989} to explore forecasting capabilities with reduced models \cite{almeida_non-intrusive_2022,Lu_databased_2017}. In contrast, our work implements a physical structure-preserving non-intrusive model learning method.
\subsubsection{Setup}\label{subsubsec:4-kse-setup}
\begin{figure}[b!]
    \centering
    \includegraphics[width=\textwidth]{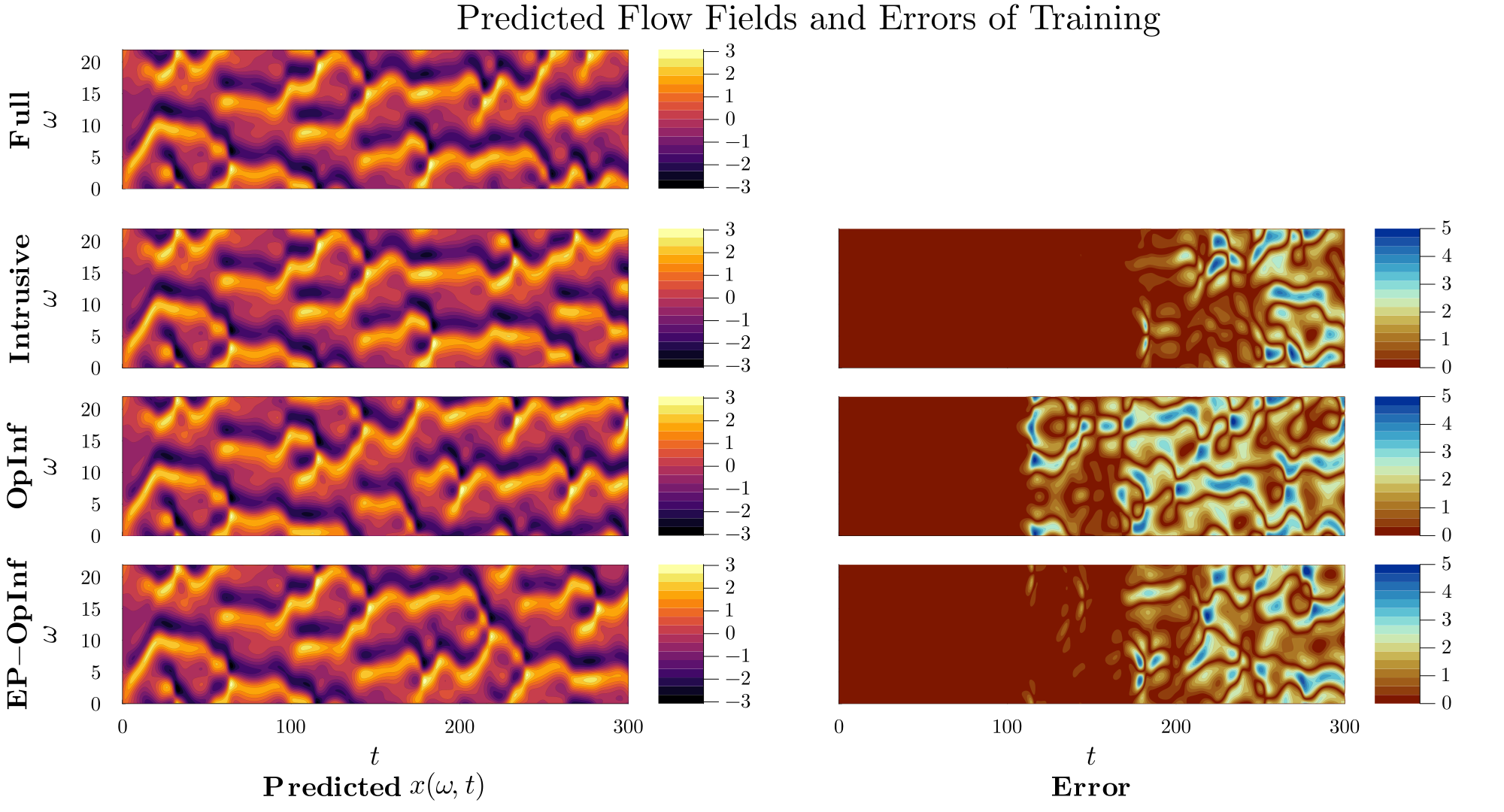}
    \caption{Training data flow field plots of the full and reduced models for an initial condition with $a=0.8,~b=0.2$ and a reduced dimension of $r=24$. The first column shows the state flow fields of the full model and reconstructed state flow fields for each method. The second column illustrates the errors of flow fields compared to the full model states.}
    \label{fig:train-kse-flow-field}
\end{figure}

We solve the KSE on the spatial domain $\omega \in [0, L)$ and the time domain $t\in[0,T]$, with periodic boundaries. 
We set $\mu = 1$, and consider the initial condition $x(\omega, 0) = a\cos(2\pi\omega/L)+b\cos(4\pi\omega/L)$, where $a$ and $b$ are parameters. 
For our experiment, we chose the spatial domain size of $L=22$, which, according to \cite{cvitanovic_state_2010}, is a sufficiently large domain size for the KSE's dynamics to transition to an inertial manifold characterized as a chaotic attractor \cite{hyman_order_1986}. Considering the increase in error for the intrusive model beyond $t>150$, as illustrated in Figure \ref{fig:train-kse-flow-field}, we chose $T=300$ for the time domain to adequately capture the progression of dynamics from the transient phase into the chaotic regime. 

Furthermore, we discretize this PDE on a uniform spatial grid of $n=512$ grid points, yielding an ODE in the form of \eqref{eqn:quadratic-system}. To integrate the ODE, we implement a semi-implicit scheme where we use the Crank-Nicolson scheme for the linear term and the Adams-Bashforth scheme for the quadratic term, maintaining a constant time step size of $\Delta t = 10^{-3}$. We collect training data at $a= [0.8, 1.0, 1.2]$ and $b=[0.2, 0.4, 0.6]$ for a total of 9 initial conditions, and retain the state snapshot at every 100\textsuperscript{th} timestep in the training data set.

\begin{figure}[h!]
    \centering
    \includegraphics[width=0.47\textwidth]{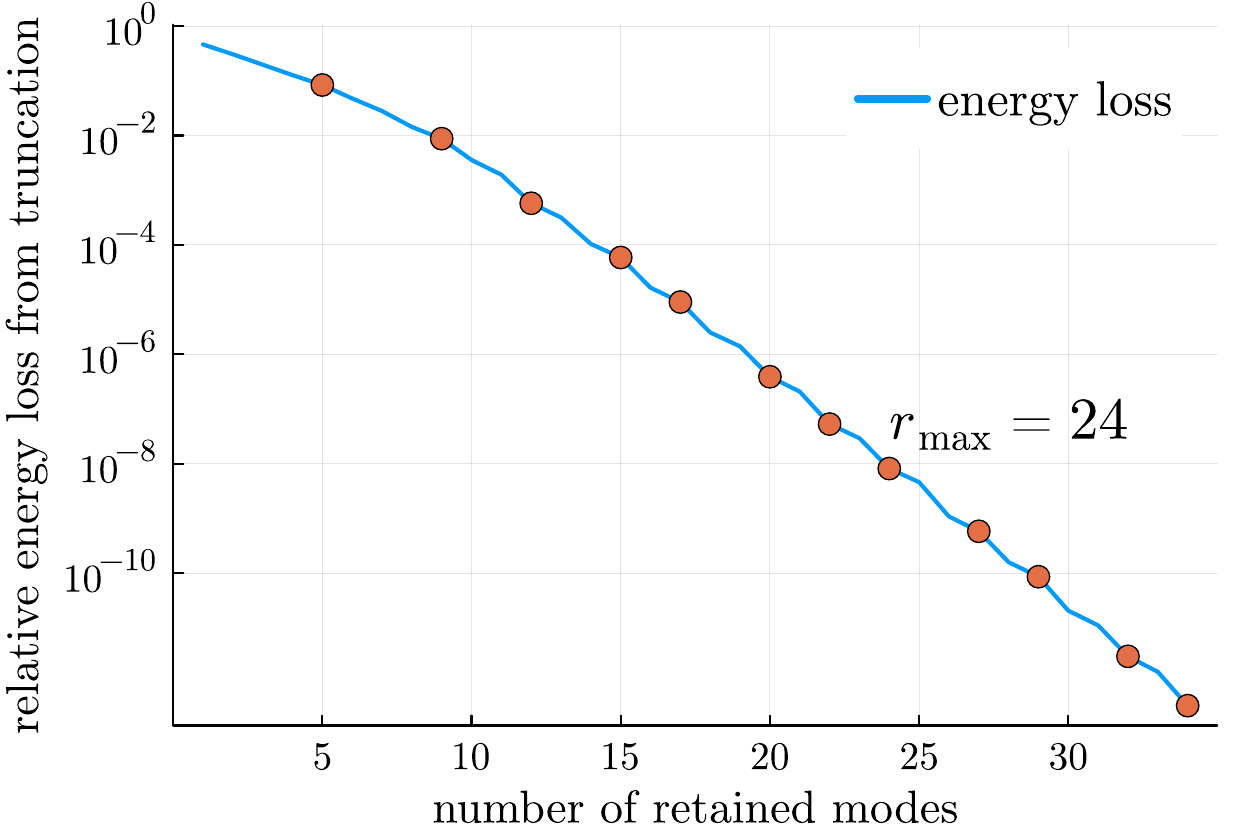}
    \caption{Relative energy lost from truncated modes for Kuramoto-Sivashinksy equation training data.}
    \label{fig:kse-eng-spectrum}
\end{figure}

Figure \ref{fig:kse-eng-spectrum} plots the energy spectrum~\eqref{eqn:energy-spectrum} of the KSE data set. We note that $r = 9$ is required to retain 99\% of the snapshot energy, while $r = 24$ retains 99.999999\%. Therefore, for both standard OpInf and EP-OpInf, the models are learned using $r_{\text{max}}=24$, and we extract submatrices of sizes $r>9$ and $r<r_{max}$. We note that the POD basis sizes required to retain most of the data energy are notably higher than those for the Burgers' equation. This is expected due to the chaotic/turbulent nature of the KSE.

After learning the models, we test them on two different test data sets not seen in training:
\begin{itemize}
    \item \textbf{Test 1 (Interpolation)}: 50 uniformly sampled initial conditions within the training parameter region: $a\in [0.8, 1.2]$ and $b\in[0.2, 0.6]$
    \item \textbf{Test 2 (Extrapolation)}: 50 uniformly sampled initial conditions outside the training parameter region: $a\in[0, 0.8) \cup (1.2, 2.0]$ and $b \in [0, 0.2) \cup (0.6, 0.8]$.
\end{itemize}

\subsubsection{Measures of reduced model performance}
The chaotic nature of the KSE dynamics means that small perturbations can lead to large differences in the state at a later time. Therefore, instead of assessing the reduced model performance using pointwise error metrics, we instead consider the reduced model predictions of the system statistics.  We follow the work \cite{Lu_databased_2017} which uses the state autocorrelation function as a statistical measure for assessing the performance of surrogate models. 
The normalized autocorrelation function is defined as follows:
\begin{equation}\label{eqn:autocorrelation}
    \rho(k) = \frac{\mathrm{Cov}[x(t),x(t+k)]}{\sqrt{\mathrm{Var}[x(t)]\mathrm{Var}[x(t+k)]}} ~,
\end{equation}
where $k$ is the time lag and the covariance and variance are taken over infinite time \cite{box_timeseries_1994}.  The autocorrelation quickly decays to zero as the time lag increases, a phenomenon known as correlation splitting which indicates the statistical independence of the system state at a given time $t$ from its past and future trajectories and is a key indicator of chaos \cite{billingsley_ergodic_1965,Cornfeld1982,Zaslavskii1985-aq,sinai_stochasticity_1979,anishchenko_correlation_2003}.

In the case of finite time series data of a stationary process, we use the \textit{sample} normalized autocorrelation function 
\begin{equation}\label{eqn:sample-autocorrelation}
    \hat \rho(k) = \frac{c_k}{c_0} \quad \text{ where ~~} c_k = \frac{1}{T}\sum_{t=1}^{T-k}\left(x(t) - \bar x\right)\left(x(t+k)-\bar x\right) ~.
\end{equation}
Here, $\bar x$ is the sample mean or time-averaged state values \cite{box_timeseries_1994}.
To quantitatively evaluate our reduced models' performance, we introduce the normalized autocorrelation error, averaged over $J$ initial conditions:
\begin{equation}\label{eqn:nace}
    \text{average normalized autocorrelation error} = \frac{1}{J}\sum_{j=1}^J\frac{\|\hat\rho_j(k) - \bar{\hat\rho}_j(k)\|^2_2}{\| {\hat\rho}_j(k)\|^2_2},
\end{equation}
where $\hat\rho_j(k)$ and $\bar{\hat\rho}_j(k)$ represent the sample normalized autocorrelation of the full and reduced models, respectively. 

We thus report the autocorrelation \eqref{eqn:autocorrelation} and autocorrelation error \eqref{eqn:nace} for both the learned and intrusive models in our results. Additionally, to qualitatively assess the model predictions, we present flow field plots alongside their corresponding errors from the full model.
\begin{figure}[t!]
    \centering
    \includegraphics[width=\textwidth]{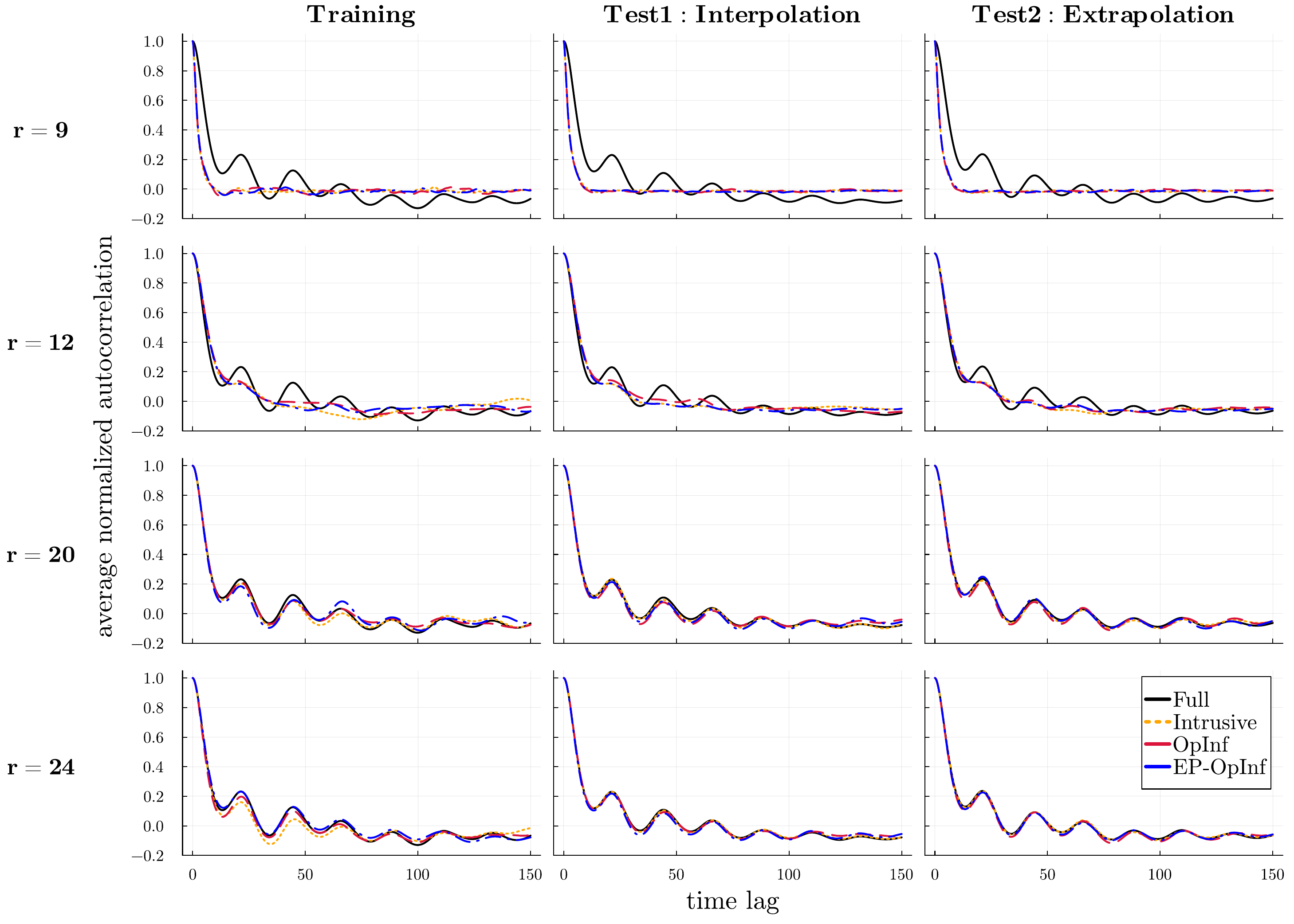}
    \caption{Normalized autocorrelation averaged over all training initial conditions. The states of each normalized autocorrelation function are reconstructed with reduced dimensions of $r=9,12,20,24$.}
    \label{fig:train-kse-ac4}
\end{figure}

\subsubsection{Results}\label{subsubsec:4-kse-results}
In the first column of Figure \ref{fig:train-kse-flow-field}, the flow fields predicted by the full model are presented alongside those reconstructed by the intrusive, standard OpInf, and EP-OpInf models with a reduced dimension of $r=24$. The second column plots the errors between the full and reduced model flow fields. We observe that all four models yield qualitatively similar predictions. During the transient time period $t\in[0,100]$, the errors are small relative to the rest of the domain, demonstrating that all three methods can accurately predict transient trajectories in the pre-chaotic regime. However, the errors begin to grow as the system transitions into the chaotic regime beyond this interval. During the interval $t\in[100,200]$, the EP-OpInf method exhibits smaller errors than the standard OpInf, and its performance is on par with the intrusive methods.
Figure \ref{fig:train-kse-ac4} illustrates the average normalized autocorrelation as a function of time lag, each at different reduced model dimensions  $r=9,12,20,24$, together with the autocorrelation of the full model. Each column represents the results for the training, Test 1 (interpolation), and Test 2 (extrapolation) data. For all three results, we observe that as the model dimension $r$ increases, the reduced models converge toward the full model's autocorrelation pattern and better capture the correlation splitting. The OpInf method generally aligns with the full model's trend but exhibits some discrepancies, particularly at lower dimensions, while the EP-OpInf is closer to the intrusive method. For Tests 1 and 2, we see that all three methods have smaller discrepancies among one another compared to the training data. This observation could be attributed to the increased number (50) of trajectories for which the error is averaged in the test cases in contrast to the 9 trajectories available for training.

\begin{figure}[t!]
    \centering
    \includegraphics[width=0.5\textwidth]{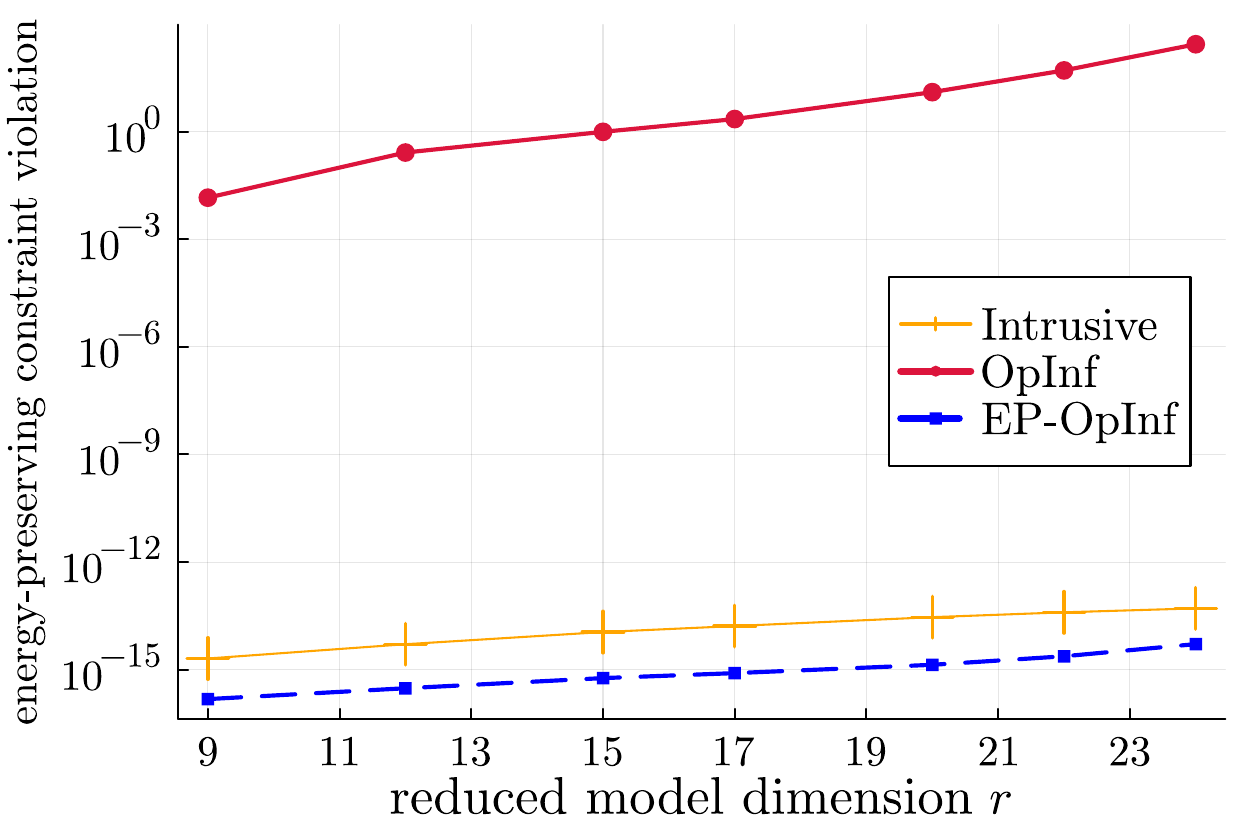}
    \caption{Energy-preserving constraint violation \eqref{eqn:constraint-violation} of the KSE.}
    \label{fig:train-kse-cv}
\end{figure}

\begin{figure}[h!]
    \centering
    \includegraphics[width=\textwidth]{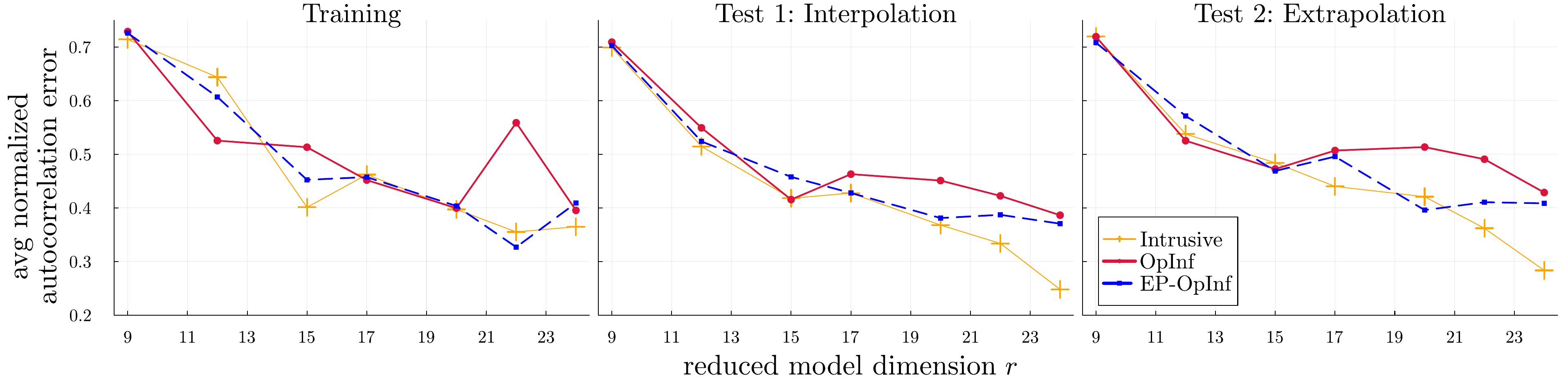}
    \caption{Normalized autocorrelation error results of the training and two tests. Training and tests are averaged over 9 and 50 initial conditions, respectively.}
    \label{fig:train-test-kse-nace}
\end{figure}

The energy-preserving constraint violation plot in Figure \ref{fig:train-kse-cv} validates that EP-OpInf allows the quadratic operator to satisfy the energy-preserving structure. Furthermore, we observe a decreasing trend in the autocorrelation error plots of Figure \ref{fig:train-test-kse-nace}, indicating the autocorrelation functions of the reduced models converging to that of the full model. Despite the three methods showing comparable error results, the error margins for Test 2 are slightly larger than those for Test 1 since the initial conditions are sampled outside the training parameter regimes. 

\begin{figure}[b!]
    \centering
    \includegraphics[width=\textwidth]{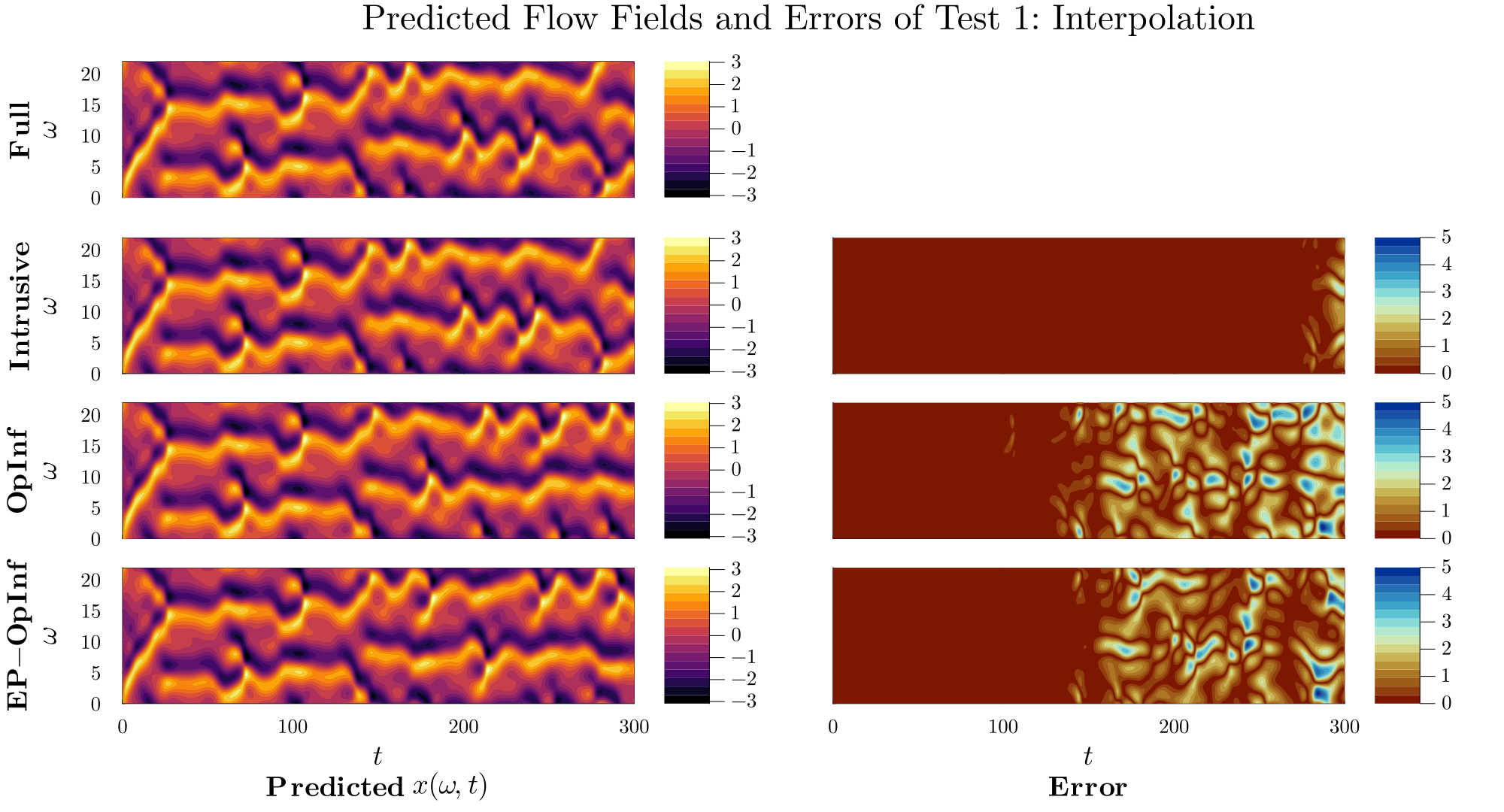}\vspace{1.5cm}
    \includegraphics[width=\textwidth]{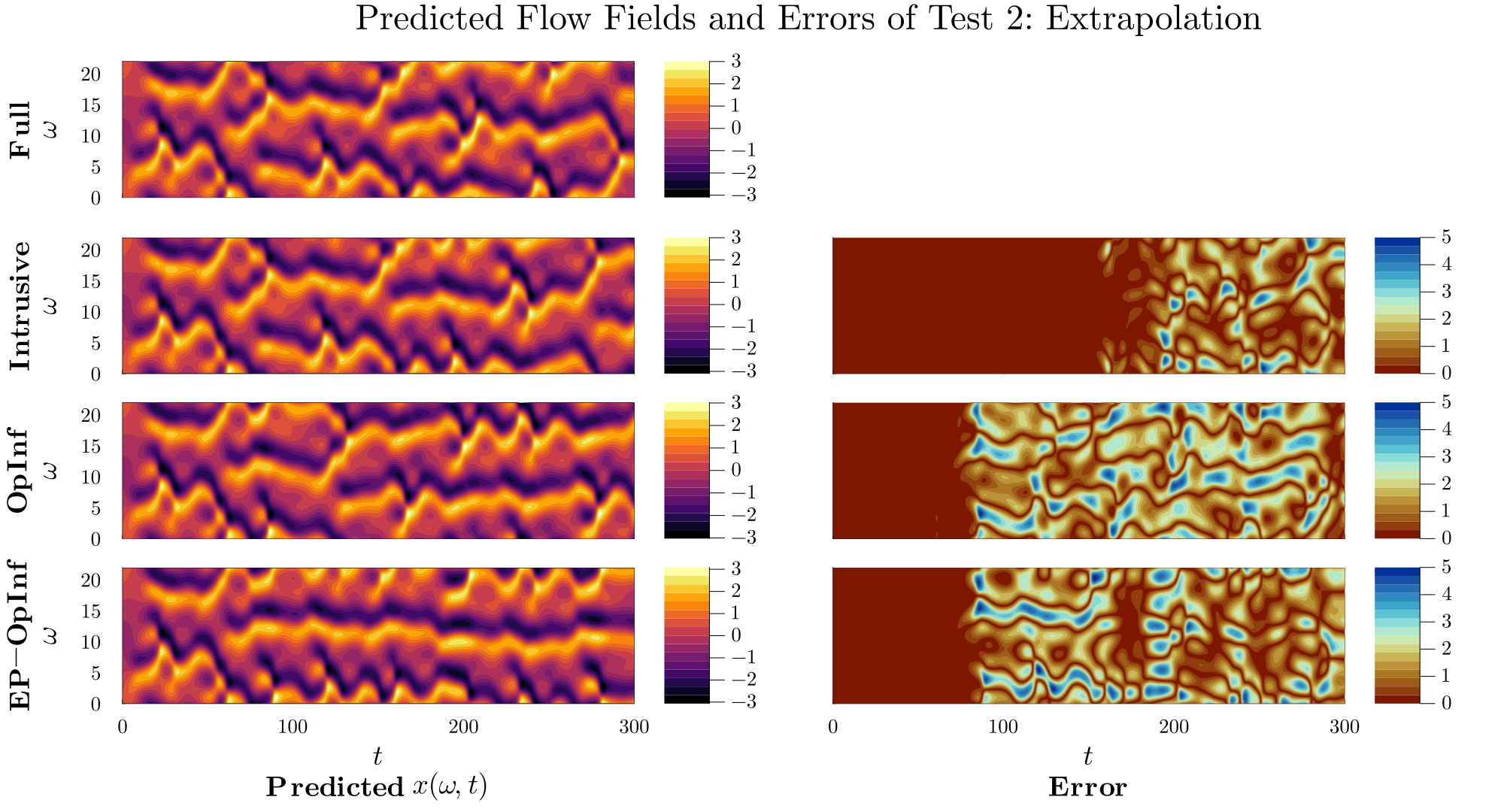}
    \caption{Flow field plots generated using one of the uniformly generated initial condition parameters for each Test 1 and Test 2. The parameters are (Test 1) $a=1.1141,~b=0.5135$ and (Test 2) $a=0.5055,~b=0.0308$, and the reduced dimension is $r=24$. The first column shows the flow fields of the full model and reconstructed state flow fields for each reduced-order method. The second column illustrates the errors of the flow fields compared to the full model states.}
    \label{fig:test-kse-flowfield}
\end{figure}

Finally, Figure \ref{fig:test-kse-flowfield} provides a visual comparison of flow fields and their corresponding errors for our two distinct tests. In Test 1, where the initial conditions are uniformly drawn from the same parameter space as the training set, the state reconstructions by EP-OpInf are on par with those by the standard OpInf method. The intrusive method, on the other hand, demonstrates very small errors. Notably, the errors in this test for standard OpInf and EP-OpInf do not show a large discrepancy from the training results demonstrated in Figure \ref{fig:train-kse-flow-field}, indicating consistent model performance within the training parameter space. In contrast, for Test 2, which involves parameters outside the training regime, OpInf and EP-OpInf exhibit larger errors across an expanded domain, but the errors of the two non-intrusive approaches are comparable to each other. Overall, EP-OpInf results demonstrate our proposed method's capability of capturing the statistical properties and model accuracy comparable to the standard OpInf and intrusive approaches while satisfying the energy-preserving constraint and being a non-intrusive method. 


\section{Conclusion}\label{sec:5-conclusion}
In this paper, we introduced Energy-Preserving Operator Inference, a new reduced model learning method for quadratic dynamical systems whose quadratic term preserves energy. Such systems arise naturally in fluid dynamics and this structure enables certain types of stability guarantees. Learning models from data which preserve this structure is therefore desirable when the learned models will be used for design and control.  The EP-OpInf method incorporates an equality constraint \eqref{eqn:energy-preservation-H} into the model learning problem to learn a reduced quadratic operator which has energy-preserving structure.

We conducted numerical experiments on model problems drawn from fluid dynamics, namely the viscous Burgers' and Kuramoto-Sivashinksy (KSE) equations. In these experiments, EP-OpInf achieved similar accuracy to the standard OpInf approach as well as to the intrusive reduced model. We emphasize that the EP-OpInf and standard OpInf approaches do this non-intrusively, which is an advantage in settings where intrusive access to the model code is not possible. Additionally, while the standard OpInf method fails to learn energy-preserving quadratic operators, our new EP-OpInf method succeeds at non-intrusively enforcing this structure. 
Directions for future work include enforcing other types of structure within the OpInf framework as well as building on the current EP-OpInf method to develop stability guarantees for learned models with energy-preserving quadratic terms. 


\bibliography{refs}

\end{document}